\documentclass[11pt,draft]{article}
\usepackage{amsmath,amscd}
\usepackage{amssymb,latexsym,amsthm}
\usepackage{color}
\usepackage[spanish,english]{babel}
\usepackage{amssymb}
\usepackage{color}
\usepackage{amsmath,amsthm,amscd}
\usepackage[latin1]{inputenc}
\usepackage{lscape}
\usepackage{fancyhdr}
\usepackage{amsfonts}
\usepackage{pb-diagram}

\numberwithin{equation}{section}
\newtheorem{theorem}{Theorem}[section]

\newtheorem{proposition}[theorem]{Proposition}

\newtheorem{corollary}[theorem]{Corollary}

\theoremstyle{definition}
\newtheorem{example}[theorem]{Example}
\newtheorem{definition}[theorem]{Definition}
\newtheorem{examples}[theorem]{Examples}

\newtheorem{remark}[theorem]{Remark}

\newcommand{\cU}{\mbox{${\cal U}$}}

\newcommand{\cW}{\mbox{${\cal W}$}}

\title{\textbf{Burchnall-Chaundy theory for skew Poincar\'e-Birkhoff-Witt extensions}}
\author{Armando Reyes\footnote{Department of Mathematics. e-mail: mareyesv@unal.edu.co} \\ Universidad Nacional de Colombia \\ H\'ector Su\'arez\footnote{School of Mathematics and Statistics. e-mail: hector.suarez@uptc.edu.co} \\ Universidad Pedag\'ogica y Tecnol\'ogica de Colombia}
\date{}

\begin{document}
\maketitle
\begin{abstract}
\noindent In this paper we review some classical results on the algebraic dependence of commuting elements in several noncommutative algebras as differential operator rings and Ore extensions. Then we extend all these results to a more general setting, the family of noncommutative rings known as skew Poincar\'e-Birkhoff-Witt extensions. 

\bigskip


\bigskip

\bigskip

\end{abstract}

\section{Introduction}
The algebraic dependence of commuting elements in noncommutative rings has its origin in a serie of papers by the British mathematicians Joseph Burchnall and Theodore Chaundy \cite{BurchnallChaundy1923}, \cite{BurchnallChaundy1928} and \cite{BurchnallChaundy1931} about commutative ordinary differential operators. Briefly, the idea is the following: Consider a commutative ring $k$ and an $k$-algebra, say $R$. Suppose that there are two commuting elements of $R$, say $a,b$. They were interested in the question whether they are algebraically dependent over $k$. More exactly, they asked if there exists  a non-zero polynomial $f(s,t)\in k[s,t]$ such that $f(a,b)=0$. \\

Several authors have considered the above question for some examples of algebras. For instance, Amitsur \cite{Amitsur1958} established results for differential operator rings on fields. A remarkable fact is that Goodearl \cite{Goodearl1983} extended the Amitsur's results to a differential operator rings  on semiprime rings. Similarly, Carlson and Goodearl \cite{CarlsonGoodearl1980} proved similar results obtained by Amitsur and Goodearl in a different setting. With respect to Ore extensions, Silvestrov et. al., \cite{Silvestrovetal2009}, \cite{HellstromSilvestrov2000}, \cite{LarssonSilvestrov2003}  have extended the result of Burchnall and Chaundy to the $q$-Weyl algebra. In fact, Mazorchuk \cite{Mazorchuk2001} has presented an alternative approach for the algebraic dependence of commuting elements in $q$-Weyl algebras. Finally, Richter \cite{Richter2014}, \cite{RichterSilvestrov2012} have established other results to Ore extensions. These results concerning about an algorithmic method to determine the polynomial $f(s,t)$ mentioned above.\\

The noncommutative rings of interest for us in this article are the skew Poincar\'e-Birkhoff-Witt extensions introduced in \cite{LezamaGallego2011}. In the literature (c.f. \cite{Reyes2013PhD} and \cite{LezamaReyes2014}) it has been showed that skew PBW extensions also generalize several families of noncommutative rings of interest in representation theory, noncommutative algebraic geometry and mathematical physics, such as the following (see \cite{ReyesSuarez2017Clifford} or \cite{ReyesSuarez2018Beit} for a detailed reference of every algebra): iterated Ore extensions of injective type defined by Ore, almost normalizing extensions defined by McConnell and Robson, solvable polynomial rings introduced by Kandri-Rody and Weispfenning, diffusion algebras, the skew polynomial rings studied by Kirkman, Kuzmanovich, and Zhang, 3-dimensional skew polynomial algebras \cite{ReyesSuarez2017PBWBasesFarEast}, and others. The importance of skew PBW extensions is that they do not assume that the coefficients commute with the variables and take coefficients not necessarily in fields. In fact, skew PBW extensions contain well-known groups of algebras such as some types Auslander-Gorenstein rings, some Calabi-Yau and skew Calabi-Yau algebras, some Artin-Schelter regular algebras, some Koszul algebras, quantum polynomials, some quantum universal enveloping algebras, some examples of $G$-algebras (see \cite{ReyesSuarez2017Clifford}, \cite{ReyesSuarez2017CalabiFarEast} and \cite{SuarezReyes2017KoszulFarEast} for more details) and many other algebras of interest for modern mathematical physicists (see \cite{Reyes2013PhD}, 
\cite{LezamaReyes2014} and \cite{ReyesSuarez2017Momento} for a list of examples). For several relations between skew PBW extensions and another algebras with PBW bases, see \cite{Reyes2013PhD}, \cite{LezamaReyes2014}, \cite{ReyesSuarez2017PBWBasesFarEast}, and \cite{SuarezReyes2017KoszulFarEast}.\\

Having in mind the general aspect of skew PBW extensions, the purpose of this paper is to establish  a first approach to the Burchnall-Chaundy theory for the skew Poincar\'e-Birkhoff-Witt extensions. With this in mind, the paper is organized as follows: we start by introducing the definition and some facts of these extensions. Then we proceed to review results obtained by other authors for particular examples of skew PBW extensions mentioned above (differential operator rings, $q$-Weyl algebras, and Ore extensions), and finally we describe the Burchnall-Chaundy theory for these noncommutative rings generalizing the results presented in \cite{Richter2014}. We consider that this study enriches the investigation of ring and homological properties of these extensions in the sense of \cite{Reyes2013PhD}, \cite{LezamaReyes2014},  \cite{ReyesSuarez2017PBWBasesFarEast}, \cite{ReyesSuarez2017Clifford}, 
\cite{ReyesSuarez2017CalabiFarEast},  \cite{SuarezReyes2017KoszulFarEast}, 
\cite{SuarezReyes2017JP}, 
\cite{Reyes2018} and \cite{ReyesSuarez2018Beit}.\\

Throughout this paper the rings and algebras are associative with unit and all modules are unital right modules. All morphisms between rings are assumed to map the multiplicative identity element to the multiplicative identity element. If $R$ is a ring and $a$ is an element of $R$, the {\em centralizer} of $a$, denoted $C_R(a)$, is the set of all elements in $R$ that commute with $a$. Now, if $k$ is  a commutative ring and $R$ is an $k$-algebra, two commuting elements $p, q\in R$ are said to be \textit{algebraically dependent} (over $k$), if there is a non-zero polynomial $f(s,t)\in k[s,t]$ such that $f(p,q)=0$, in which case $f$ is called an annihilating polynomial.  
\section{Definitions and elementary properties}\label{definitionexamplesspbw}
We start recalling the definition and some properties about our object of study.
\begin{definition}[\cite{LezamaGallego2011}, Definition 1]\label{gpbwextension}
Let $R$ and $A$ be rings. We say that $A$ is a \textit{skew Poincar\'e-Birkhoff-Witt extension of} $R$, briefly skew PBW extension, if the following conditions hold:
\begin{enumerate}
\item[\rm (i)]$R\subseteq A$;
\item[\rm (ii)]there exist elements $x_1,\dots ,x_n\in A\ \backslash\ R$ such that $A$ is a left free $R$-module, with basis the basic elements ${\rm Mon}(A):= \{x^{\alpha}=x_1^{\alpha_1}\cdots
x_n^{\alpha_n}\mid \alpha=(\alpha_1,\dots ,\alpha_n)\in \mathbb{N}^n\}\ (x^{0}_1\dotsb x_n^{0}:=1)$.
\item[\rm (iii)]For each $1\leq i\leq n$ and any $r\in R\ \backslash\ \{0\}$, there exists an element $c_{i,r}\in R\ \backslash\ \{0\}$ such that $x_ir-c_{i,r}x_i\in R$.
\item[\rm (iv)]For any elements $1\leq i,j\leq n$,  there exists $c_{i,j}\in R\ \backslash\ \{0\}$ such that $x_jx_i-c_{i,j}x_ix_j\in R+Rx_1+\cdots +Rx_n$. 
\end{enumerate}
Under these conditions we write $A:=\sigma(R)\langle
x_1,\dots,x_n\rangle$.
\end{definition}
\begin{proposition}[\cite{LezamaGallego2011}, Proposition
3]\label{sigmadefinition}
Let $A$ be a skew PBW extension of $R$. For each $1\leq i\leq
n$, there exist an injective endomorphism $\sigma_i:R\rightarrow
R$ and a $\sigma_i$-derivation $\delta_i:R\rightarrow R$ such that $x_ir=\sigma_i(r)x_i+\delta_i(r)$, for any $r \in R$.
\end{proposition}
\begin{definition}\label{setsinjectivederivations}
Let $\Sigma:=\{\sigma_1,\dotsc, \sigma_n\}$ and $\Delta:=\{\delta_1,\dotsc, \delta_n\}$. Any element $r$ of $R$ such that $\sigma_i(r)=r$ and $\delta(r)=0$ will be called a {\em constant}. Note that the constants form a subring of $R$, This ring will be denoted by $F$.
\end{definition}
\begin{example}\label{mentioned}
Remarkable examples of skew PBW extensions are the differential operator rings and the Ore extensions of injective type, see \cite{LezamaReyes2014} and  \cite{Reyes2013PhD} for more details. These examples include algebras of interest for modern mathematical physicists such  group rings of polycyclic-by-finite groups, Ore algebras, operator algebras, diffusion algebras, some quantum groups, quadratic algebras in three  variables, and $3$-dimensional skew polynomial rings. 
\end{example}
\begin{definition}[\cite{LezamaGallego2011}, Definition 6]\label{definitioncoefficients}
Let $A$ be a skew PBW extension of $R$ with endomorphisms
$\sigma_i$, $1\leq i\leq n$, as in Proposition
\ref{sigmadefinition}.
\begin{enumerate}
\item[\rm (a)]For $\alpha=(\alpha_1,\dots,\alpha_n)\in \mathbb{N}^n$,
$\sigma^{\alpha}:=\sigma_1^{\alpha_1}\cdots \sigma_n^{\alpha_n}$,
$|\alpha|:=\alpha_1+\cdots+\alpha_n$. If
$\beta=(\beta_1,\dots,\beta_n)\in \mathbb{N}^n$, then
$\alpha+\beta:=(\alpha_1+\beta_1,\dots,\alpha_n+\beta_n)$.
\item[\rm (b)]For $X=x^{\alpha}\in {\rm Mon}(A)$,
$\exp(X):=\alpha$ and $\deg(X):=|\alpha|$. The symbol $\succeq$ will denote a total order defined on ${\rm Mon}(A)$ (a total order on $\mathbb{N}_0^n$). For an
 element $x^{\alpha}\in {\rm Mon}(A)$, ${\rm exp}(x^{\alpha}):=\alpha\in \mathbb{N}_0^n$.  If
$x^{\alpha}\succeq x^{\beta}$ but $x^{\alpha}\neq x^{\beta}$, we
write $x^{\alpha}\succ x^{\beta}$. If $f=c_1X_1+\dotsb +c_tX_t\in
A$, $c_i\in R\ \backslash\ \{0\}$, with $X_1\succ \dotsb \succ X_t$, then ${\rm
lm(f)}:=X_1$ is the \textit{leading monomial} of $f$, ${\rm
lc}(f):=c_1$ is the \textit{leading coefficient} of $f$, ${\rm
lt}(f):=c_1X_1$ is the \textit{leading term} of $f$,  ${\rm exp}(f):={\rm exp}(X_1)$ is the \textit{order} of $f$, and
 $E(f):=\{{\rm exp}(X_i)\mid 1\le i\le t\}$. Finally, if $f=0$, then
${\rm lm}(0):=0$, ${\rm lc}(0):=0$, ${\rm lt}(0):=0$. We also
consider $X\succ 0$ for any $X\in {\rm Mon}(A)$. For a detailed description of monomial orders in skew PBW extensions, see \cite{LezamaGallego2011}, Section 3. For an element $f\in A$ as above, $\deg(f):=\max\{\deg(X_i)\}_{i=1}^t$.
\end{enumerate}
\end{definition}
\begin{proposition}[\cite{LezamaReyes2014}, Proposition 4.1]\label{extensiondomain}
Let $R$ be a domain. If $A$ is a skew PBW extension of $R$, then $A$ is also a domain.
\begin{proof}
Let $f=cx^{\alpha}+p, g=dx^{\beta} + q$ be nonzero elements of $A$, with ${\rm lt}(f)=cx^{\alpha}, {\rm lt}(g)=dx^{\beta}$, so $c,d\neq 0$, and $x^{\alpha}\succ {\rm lm}(p)$ and $x^{\beta}\succ {\rm lm}(q)$. We have
\begin{align*}
fg = &\ (cx^{\alpha}+p)(dx^{\beta}+q) = cx^{\alpha}dx^{\beta} + cx^{\alpha}q +pdx^{\beta} + pq\\
= &\ c(d_{\alpha}x^{\alpha} + p_{\alpha,d})x^{\beta} + cx^{\alpha}q + pdx^{\beta} + pq,
\end{align*}
with $0\neq d_{\alpha}=\sigma^{\alpha}(d)\in R$, $p_{\alpha, d}\in A, p_{\alpha, d}=0$ or ${\rm deg}(p_{\alpha, d}) < |\alpha|$ if $p_{\alpha, d} \neq 0$. In this way
\begin{align*}
fg = &\ cd_{\alpha}x^{\alpha}x^{\beta} + cp_{\alpha, d}x^{\beta} + cx^{\alpha}q + pdx^{\beta} + pq\\
= &\ cd_{\alpha}(c_{\alpha, \beta}x^{\alpha+\beta} + p_{\alpha, \beta}) + cp_{\alpha,d}x^{\beta} + cx^{\alpha}q + pdx^{\beta} + pq\\
= &\ cd_{\alpha}c_{\alpha, \beta}x^{\alpha+\beta} + cd_{\alpha}p_{\alpha, \beta} + cp_{\alpha, d}x^{\beta} + cx^{\alpha}q + pdx^{\beta} + pq,
\end{align*}
where $0\neq c_{\alpha, \beta}\in R,\ p_{\alpha, \beta}\in A,\ p_{\alpha, \beta}=0$ or ${\rm deg}(p_{\alpha, \beta}) < |\alpha + \beta|$ if $p_{\alpha, \beta}\neq 0$. Moreover $cd_{\alpha}c_{\alpha, \beta}\neq 0$ and $h:= cd_{\alpha}p_{\alpha, \beta} + cp_{\alpha, d}x^{\beta} + cx^{\alpha}q + pdx^{\beta} + pq\in A$ is such that $h=0$ or $x^{\alpha +\beta} \succ {\rm lm}(h)$, that is, $fg\neq 0$. 
\end{proof}
\end{proposition}

\section{Burchnall-Chaundy theory for remarkable examples of skew Poincar\'e-Birkhoff-Witt extensions}
In this section we present a review of the most important results about Burchnall-Chaundy theory for remarkable examples of noncommutative rings, the differential operator rings and Ore extensions. We consider a similar presentation to the established in \cite{Richter2014}. 
\subsection{Burchnall-Chaundy theory for differential operator rings}
We start describing some results on the algebraic dependence of commuting elements in differential operator rings which are particular examples of skew PBW extensions. A proof of this assertion can be found in \cite{LezamaGallego2011},  \cite{LezamaReyes2014} or \cite{ReyesSuarez2017Momento}. Now, as we said in the Introduction, this sort of question has its origin in a series of papers by the British mathematicians Joseph Burchnall and Theodore Chaundy, see \cite{BurchnallChaundy1923}, \cite{BurchnallChaundy1928}, and \cite{BurchnallChaundy1931}.\\

Amitsur \cite{Amitsur1958} studied the case when $\Bbbk$ is a field of characteristic zero and $\delta$ is an arbitrary derivation on $\Bbbk$. He obtained the following two results.
\begin{proposition}[\cite{Amitsur1958}, Theorem 1]
Let $\Bbbk$ be a field of characteristic zero with a derivation $\delta$. Let $F$ denote the subfield of constants. Let $S=\Bbbk[x;{\rm id}_{\Bbbk}, \delta]$ be the differential operator ring, and let $P$ be an element of $S$ of degree $n$. Denote by $F[P]$ the ring of polynomials in $P$ with constant coefficients, that is, $F[P]=\{\sum_{j=0}^m b_jP^{j}\mid b_j\in F\}$. Then $C_S(P)$ is a commutative subring of $S$ and a free $F[P]$-module of rank at most $n$.
\end{proposition}
\begin{proposition}[\cite{Amitsur1958}, Corollary 2]
Let $P$ and $Q$ be two commuting elements of the skew polynomial ring $\Bbbk[x;{\rm id}_{\Bbbk}, \delta]$, where $\Bbbk$ is a field of characteristic zero. Then there is a nonzero polynomial $f(s,t)$, with coefficients in $F$, such that $f(P,Q)=0$. 
\end{proposition}
Goodearl \cite{Goodearl1983} has extended the results of Amitsur to a more general setting. We recall that a commutative ring is semiprime if and only if it has no nonzero nilpotent elements.
\begin{proposition}[\cite{Goodearl1983}, Theorem 1.2]\label{Theorem2.3paper}
Let $R$ be a semiprime commutative ring with derivation $\delta$ and assume that its ring of constants is a field, $F$. If $P$ is an operator in $S=R[x;{\rm id}_{R}, \delta]$ of positive degree $n$, where $n$ is invertible in $F$, and has an invertible leading coefficient, then $C_S(P)$ is a free $F[P]$-module of rank at most $n$.
\end{proposition}
Goodearl \cite{Goodearl1983} showed that if $R$ is a semiprime ring of positive characteristic such that the ring of constants is a field, then $R$ must be a field. In this case he proves the following result.
\begin{proposition}[\cite{Goodearl1983}, Theorem 1.11]\label{Theorem2.4paper}
Let $\Bbbk$ be a field with a derivation $\delta$, and let $F$ be its subfield of constants. If $P$ is an element of $S=\Bbbk[x;{\rm id}_{\Bbbk}, \delta]$ of positive degree $n$ and with invertible leading coefficient, then $C_S(P)$ is a free $F[P]$-module of rank at most $n^2$.
\end{proposition}
From Propositions \ref{Theorem2.3paper} and  \ref{Theorem2.4paper} we get the following result.
\begin{corollary}[\cite{Goodearl1983}, Theorem 1.13]
Let $P$ and $Q$ be commuting elements of $R[x;{\rm id}_{R}, \delta]$, where $R$ is a semiprime commutative ring, with a derivation $\delta$ such that the subring of constants is a field. Suppose that the leading coefficient of $P$ is invertible. Then there exists a non-zero polynomial $f(s,t)\in F[s,t]$ such that $f(P,Q)=0$.
\end{corollary}
Carlson and Goodearl \cite{CarlsonGoodearl1980} proved similar results to Propositions \ref{Theorem2.3paper} and \ref{Theorem2.4paper}. More exactly:
\begin{proposition}[\cite{CarlsonGoodearl1980}, Theorem 1]
Let $R$ be a commutative ring, with a derivation $\delta$ such that the ring of constants is a field $F$ of characteristic zero. Assume that, for all $a\in R$, if the set $\{b\in R\mid \delta(b)=ab\}$ contains a nonzero element, then it contains an invertible element. Let $P$ be an element of $S=R[x;{\rm id}_R, \delta]$ of positive degree $n$ with invertible leading coefficient. Then $C_S(P)$ is a free $F[P]$-module of rank at most $n$. As before, this implies that if $Q$ commutes with $P$, there exists a nonzero polynomial $f(s,t)\in F[s,t]$  such that $f(P,Q)=0$.
\end{proposition}
\subsection{Burchnall-Chaundy theory for Ore extensions}
Let $\Bbbk$ be a field and $q$ a nonzero element of $\Bbbk$ which is not a root of unity. Set $R=\Bbbk[y]$, a polynomial ring in one variable over $\Bbbk$. There is an endomorphism $\sigma$ of $R$ such that $\sigma(y)=qy$ and $\sigma(r)=r$, for all $r\in \Bbbk$. For this $\sigma$ there exists a unique $\sigma$-derivation $\sigma$-derivation $\delta$ such that $\delta(y)=1$ and $\delta(a)=0$, for all $a\in \Bbbk$. The Ore extension $R[x;\sigma,\delta]$ for this choice of $R, \sigma$ and $\delta$ is known as the first $q$-Weyl algebra (see \cite{LezamaReyes2014} or \cite{Reyes2013PhD} for more details of this algebra). Again, this algebra is a particular example of skew PBW extensions, see \cite{LezamaGallego2011},  \cite{LezamaReyes2014} or \cite{ReyesSuarez2017Momento} for a proof of this assertion. \\

Silvestrov et. al., \cite{Silvestrovetal2009}, \cite{HellstromSilvestrov2000}, \cite{LarssonSilvestrov2003}  have extended the result of Burchnall and Chaundy to the $q$-Weyl algebra. These references contain two different proofs of the fact any pair of commuting elements of $B[x;\sigma, \delta]$ are algebraically independent over $\Bbbk$. In fact, in \cite{Silvestrovetal2009} an algorithm to compute an annihilating polynomial explicitly is presented. On the other hand, Mazorchuk \cite{Mazorchuk2001} has presented an alternative approach to showing the algebraic dependence of commuting elements in $q$-Weyl algebras. Let us see it.
\begin{proposition}\label{Theorem3.1paper}
Let $\Bbbk$ be a field and $q$ an element of $\Bbbk$. Set $B=\Bbbk[y]$ and suppose that $\sum_{i=0}^{N} q^{i}\neq 0$ for any natural number $N$. Let $P$ be an element of $S=B[x;\sigma,\delta]$ of degree at least $1$. Then $C_S(P)$ is a free $\Bbbk[P]$-module of finite rank.
\end{proposition}
If $P$ is as in Proposition \ref{Theorem3.1paper} and $Q$ is any element of $B[x;\sigma, \delta]$  that commutes with $P$, then there is an annihilating polynomial $f(s,t)$ with coefficients in $\Bbbk$. The methods used to obtain Proposition \ref{Theorem3.1paper} have been generalized by Hellstrom and Silvestrov in \cite{HellstromSilvestrov2007}.
Under certain assumptions on the endomorphism $\sigma$ we have the following proposition.
\begin{proposition}[\cite{Richter2014}, Theorem 3.3]\label{Theorem3.3paper}
Let $R$ be an integral domain, $\sigma$ an injective endomorphism of $R$ and $\delta$ a $\sigma$-derivation on $B$. Suppose that the ring of constants $F$ is a field. Let $a$ be an element of $S=R[x;\sigma,\delta]$ of degree $n$, and assume that if $b$ and $c$ are two elements in $C_S(a)$ such that ${\rm deg}(b)={\rm deg}(c)=m$, then $b_m=p c_m$, where $b_m$ and $c_m$ are the leading coefficients of $b$ and $c$, respectively, and $p$ is some constant. Then $C_S(a)$ is a free $F[a]$-module of rank at most $n$.
\end{proposition}
With the purpose of illustrate  Proposition \ref{Theorem3.3paper}, Richter \cite{Richter2014} have obtained the following result.
\begin{proposition}[\cite{Richter2014}, Proposition 3.1]\label{Proposition3.1paper}
Let $\Bbbk$ be a field and set $R=\Bbbk[y]$. Let $\sigma$ be an endomorphism of $R$ such that $\sigma(r)=r$, for all $r\in \Bbbk$ and $\sigma(y)=py$, where $p(y)$ is a polynomial of degree {\rm (}in $y${\rm )} greater than $1$. Let $\delta$ be a $\sigma$-derivation such that $\delta(a)=0$ for all $a\in \Bbbk$. Consider the Ore extension $S=R[x;\sigma,\delta]$ {\rm (}note that its ring of constants is $\Bbbk${\rm )}. Let $a\notin \Bbbk$ be an element of $R[x;\sigma,\delta]$. Suppose that $b,c$ are elements of $S$ such that ${\rm deg}(b)={\rm deg}(c)=m$ {\rm (}here the degree is taken with respect to $x${\rm )}, and $b,c\in C_S(a)$. Then $b_m=\alpha c_m$, where $b_m, c_m$ are the leading coefficients of $b$ and $c$, respectively, and $\alpha$ is some constant.
\end{proposition}

If $\Bbbk, \sigma, \delta, a$ as in Proposition \ref{Proposition3.1paper}, then $C_S(a)$ is a free $\Bbbk[a]$-module of finite rank (\cite{Richter2014}, Proposition 3.2). 

\begin{proposition}[\cite{Richter2014}, Theorem 3.4]\label{Theorem3.4paper}
Let $\Bbbk$ be a field,  $\sigma$  an endomorphism of $\Bbbk[y]$ such that $\sigma(y)=p(y)$, where ${\rm deg}(p)>1$, and let $\delta$ be a $\sigma$-derivation. Suppose that $\sigma(r)=r$ and $\delta(r)=0$, for all $r\in \Bbbk$. Let $a,b$ be two commuting elements of $\Bbbk[y][x;\sigma,\delta]$. Then there is a nonzero polynomial $f(s,t)\in \Bbbk[s,t]$ such that $f(a,b)=0$.
\end{proposition}

\begin{proposition}[\cite{HellstromSilvestrov2000}, Theorem 7.5]
Let $R=\Bbbk[y]$, $\sigma(y)=qy$ and $\delta(y)=1$, where $q\in \Bbbk$ and $q$ is a root of unity. Let $S=R[x;\sigma, \delta]$ and  $C$ the center of $S$. If $a,b$ are commuting elements of $S$ then there is a nonzero polynomial $f(s,t)\in C[s,t]$ such that $f(a,b)=0$.
\end{proposition}
\section{Burchnall-Chaundy theory for skew Poincar\'e-Birkhoff-Witt extensions}
In this section we present an extension of Burchnall-Chaundy theory to skew PBW extensions. Recall that $F$ is the subring of constants of $R$, and $F[P]$ is the ring of polynomials in $P$ with constant coefficients, that is, $F[P]=\{\sum_{j=0}^m b_jP^{j}\mid b_j\in F\}$.

We start with the following theorem which is the first important result of the paper.
\begin{theorem}\label{Theorem2.2}
Let $R$ be an integral domain and let $A$ be a skew PBW extension of $R$. Suppose that the ring of constants $F$ is a field. Let $f\in A$ with ${\rm exp}(f)=\alpha$, and assume that if $g, h\in C_A(f)$ with $|{\rm exp}(g)|=|{\rm exp}(h)|$, then ${\rm lc}(g)=p {\rm lc}(h)$, where $p$ is some constant. Then $C_A(f)$ is a free $F[f]$-module of rank at most $|\alpha|$.\ 
\begin{proof}
We follow similar ideas to the used in the proof or Proposition \ref{Theorem3.3paper}. Let $M$ be the subset of elements of $\{0,1,\dotsc, |\alpha|-1\}$ such that an integer $0\le i< |\alpha|$ is an element of $M$ if and only if $C_A(f)$ contains an element of degree equivalent to $i$ modulo $|\alpha|$. For $i\in M$ let $p_i$ be an element in $C_A(f)$ such that ${\rm deg}(p_i)\equiv i ({\rm mod}\ |\alpha|)$ and $p_i$ has minimal degree for this property. Let $p_0:=1$. The idea is to show that the set $\{p_i\mid i\in M\}$ is a basis for $C_A(f)$ as a $F[f]$-module.

Since $A$ is an integral domain (Proposition \ref{extensiondomain}) and every $\sigma_i$, $1\le i\le n$, is injective (Proposition \ref{sigmadefinition}), the degree of a product of two elements in $A$ is the sum of the degrees of the two elements. First of all, we prove that the $p_i$ are linearly independent over $F[f]$. Suppose $\sum_{i\in M}f_ip_i=0$, for some $f_i\in F[f]$. If $f_i\neq 0$, then ${\rm deg}(f_i)$ is divisible by $|\alpha|$, and hence
\begin{equation}\label{(2)}
{\rm deg}(f_ip_i)={\rm deg}(f_i) + {\rm deg}(p_i)\equiv {\rm deg}(p_i)\equiv i ({\rm mod}\ |\alpha|).
\end{equation}
Now, if $\sum_{i\in M} f_ip_i=0$ but not all $f_i$ are zero, then we have two nonzero terms  $f_ip_i$ and $f_jp_j$ with the same degree, even $i,j\in M$ distinct. Nevertheless, this can not be the case since the relation $i\equiv j ({\rm mod}\ |\alpha|)$ is false.

Second of all, let us see that the $p_i$'s span $C_A(f)$. Let $T$ denote the submodule they do span. The proof is by induction on the degree to show that all elements of $C_A(f)$ belong to $T$. Note that if $e$ is an element of degree $0$ in $C_A(f)$, by assumption on $f$ applied to $e$, we find that $e=t$ for some $t\in F$ (recall that $p_0=1$). Therefore $e\in T$.

Suppose that $T$ contains all elements in $C_A(f)$ of degree less than $j$, and let $e\in C_A(f)$ with ${\rm exp}(e)=\gamma$ with $|\gamma|=j$. There exists $i\in M$ with $j\equiv i ({\rm mod}\ |\alpha|)$. Set ${\rm deg}(p_i):=m$. We have $m\equiv j ({\rm mod}\ |\alpha|)$ and $m\le j$. In this way $j=m+q|\alpha|$ for some non-negative integer $q$. Note that  $f^{q}p_i\in T$  and ${\rm deg}(f^{q}p_i)=j$. By assumption, ${\rm lc}(e)=p {\rm lc}(f^qp_i)$, $p$ a constant. Hence, the element $e-pf^{q}p_i$ lies in $C_A(f)$ and ${\rm deg}(e-pf^{q}p_i) < j$. By the induction hypothesis it also lies in $T$, and therefore $e\in T$.
\end{proof}
\end{theorem}
Proposition \ref{primeracercamiento} establishes  sufficient conditions to use Theorem \ref{Theorem2.2}. 
\begin{proposition}\label{primeracercamiento}
Let $\Bbbk$ be a field and $A$ a  skew PBW extension of\ $\Bbbk$. If $f\in A\ \backslash\ \Bbbk$, $g, h\in C_A(f)$ with ${\rm exp}(g)={\rm exp}(h)=\beta$, and ${\rm lc}(g), \frac{1}{{\rm lc}(h)}$ are constants, then ${\rm lc}(g)=p\ {\rm lc}(h)$, where $p$ is a constant.
\begin{proof}
Consider $f={\rm lc}(f)x^{\alpha}+q_{\alpha}$, $g={\rm lc}(g)x^{\beta}+q_{\beta}$, where ${\rm deg}(q_{\alpha})< {\rm deg}(x^{\alpha})$, and ${\rm deg}(q_{\beta}) < {\rm deg}(x^{\beta})$. Since $fg=gf$ we have 
\begin{align*}
({\rm lc}(f)x^{\alpha}+q_{\alpha})({\rm lc}(g)x^{\beta} +q_{\beta}) = &\ ({\rm lc}(g)x^{\beta} + q_{\beta})({\rm lc}(f)x^{\alpha} + q_{\alpha})\\
{\rm lc}(f)x^{\alpha}{\rm lc}(g)x^{\beta} + {\rm other\ terms}  = &\ {\rm lc}(g)x^{\beta}{\rm lc}(f)x^{\alpha} + {\rm other\ terms}\\
{\rm lc}(f)[\sigma^{\alpha}({\rm lc}(g))x^{\alpha} + p_{\alpha, {\rm lc}(g)}]x^{\beta} + {\rm other\ terms} = &\ {\rm lc}(g)[\sigma^{\beta}({\rm lc}(f))x^{\beta}+p_{\beta, {\rm lc}(f)}]x^{\alpha} + {\rm other\ terms}
\end{align*}
where ${\rm deg}(p_{\alpha, {\rm lc}(g)})=0$, or ${\rm deg}(p_{\alpha, {\rm lc}(g)}) < {\rm deg}(x^{\alpha})$ if $p_{\alpha, {\rm lc}(g)}\neq 0$, and ${\rm deg}(p_{\beta, {\rm lc}(f)})=0$, or ${\rm deg}(p_{\beta,{\rm lc}(f)}) < {\rm deg}(x^{\beta})$ if $p_{\beta, {\rm lc}(f)}\neq 0$. Now
\begin{align*}
{\rm lc}(f)\sigma^{\alpha}({\rm lc}(g))x^{\alpha}x^{\beta} + {\rm other\ terms} = &\ {\rm lc}(g)\sigma^{\beta}({\rm lc}(f))x^{\beta}x^{\alpha} + {\rm other\ terms}\\
{\rm lc}(f)\sigma^{\alpha}({\rm lc}(g))[c_{\alpha, \beta}x^{\alpha+\beta} + p_{\alpha, \beta}] + {\rm other\ terms} = &\ {\rm lc}(g)\sigma^{\beta}({\rm lc}(f))[c_{\beta, \alpha}x^{\alpha+\beta}+p_{\beta, \alpha}] + {\rm other\ terms}
\end{align*}
with $p_{\alpha, \beta}=0$, or ${\rm deg}(p_{\alpha, \beta}) < |\alpha + \beta|$, and $p_{\beta, \alpha}=0$, or ${\rm deg}(p_{\beta, \alpha}) < |\alpha + \beta|$.  In this way
\begin{equation}\label{equation3commute}
{\rm lc}(f)\sigma^{\alpha}({\rm lc}(g))c_{\alpha, \beta} = {\rm lc}(g)\sigma^{\beta}({\rm lc}(f))c_{\beta, \alpha}.
\end{equation}
Similarly, one can show that the equality $ad=da$ imply that 
\begin{equation}\label{equation4commute}
{\rm lc}(f)\sigma^{\alpha}({\rm lc}(h))c_{\alpha, \beta} = {\rm lc}(h) \sigma^{\beta}({\rm lc}(f)) c_{\beta, \alpha}.
\end{equation}\ \ 
If we divide (\ref{equation3commute}) into (\ref{equation4commute}) we obtain
\begin{equation}\label{obligaanoserbijective}
\frac{\sigma^{\alpha}({\rm lc}(g))}{\sigma^{\alpha}({\rm lc}(h))} = \frac{{\rm lc}(g)}{{\rm lc}(h)},\ \ {\rm or\ equivalently}\ \ {\rm lc}(g)= \frac{\sigma^{\alpha}({\rm lc}(g))}{\sigma^{\alpha}({\rm lc}(h))} {\rm lc}(h)
\end{equation}
Using the fact that $\sigma^{\alpha}$ is an endomorphism, the expression (\ref{obligaanoserbijective}) can be written as ${\rm lc}(g) = \sigma^{\alpha}\biggl(\frac{{\rm lc}(g)}{{\rm lc}(h)}\biggr) {\rm lc}(h)$. By assumption, the elements ${\rm lc}(g)$ and $\frac{1}{{\rm lc}(h)}$ are constants, so $\sigma^{\alpha}\biggl(\frac{{\rm lc}(g)}{{\rm lc}(h)}\biggr) = \frac{{\rm lc}(g)}{{\rm lc}(h)}$ whence it is clear that $\sigma_i\biggl(\sigma^{\alpha}\biggl(\frac{{\rm lc}(g)}{{\rm lc}(h)}\biggr)\biggr) = \sigma^{\alpha}\biggl(\frac{{\rm lc}(g)}{{\rm lc}(h)}\biggr)$, for every $1\le i\le n$. Finally, since $\delta_i\biggl(\frac{{\rm lc}(g)}{{\rm lc}(h)}\biggr) = {\rm lc}(g) \delta_i(\frac{1}{{\rm lc}(h)}) + \delta_i({\rm lc}(g)) \frac{1}{{\rm lc}(h)} = 0$, since ${\rm lc}(g)$ and $\frac{1}{{\rm lc}(h)}$ are constants, then $\delta_i\biggl(\frac{{\rm lc}(g)}{{\rm lc}(h)}\biggr) = \delta_i\biggl(\sigma^{\alpha}\biggl(\frac{{\rm lc}(g)}{{\rm lc}(h)}\biggr)\biggr) = 0$, so the assertion follows.
\end{proof}
\end{proposition}
The following theorem is the more important result of this paper.
\begin{theorem}\label{moreimportant}
Let $P$ and $Q$ be two commuting elements of $A$, where $A$ is a skew PBW extension of a field $\Bbbk$. Then there is a polynomial $f(s,t)$ with coefficients in $F$ satisfying the condition $f(P,Q)=0$. 
\begin{proof}
Let $P$ be an element of $A$. Consider $P={\rm lc}(P)x^{\alpha}+q_P$ with ${\rm deg}(q_p) < {\rm deg}(x^{\alpha})$. By assumption $Q\in C_A(P)$, which implies that the elements $1, Q, Q^2, \dotsc, Q^{|\alpha|}$ are also elements of $C_A(P)$. Now, Theorem \ref{Theorem2.2} imply  that these elements are linearly dependent over $F[P]$. In other words, there are (not all) nonzero elements $g_0(P), g_1(P), \dotsc, g_{|\alpha|}(P)\in F[P]$ with
\[
g_0 + g_1Q + \dotsb + g_{|\alpha|} Q^{|\alpha|} = 0.
\]
Then, if one consider $f(s,t):=\sum_{i=0}^{|\alpha|} g_i(s)t^{i}$, one obtain that $f(P,Q)=0$.
\end{proof}
\end{theorem}
\begin{examples}
Theorem \ref{moreimportant} can be applied to several skew PBW extensions. More precisely, if $A$ is a skew PBW extension of a field $\Bbbk$ where the coefficients commute with the variables, that is, $x_ir = rx_i$, for every $r\in R$ and each $i=1,\dotsc, n$ (these extensions were called {\em constant} by the authors in \cite{SuarezReyes2017KoszulFarEast}), then Theorem \ref{moreimportant} can be illustrated. Some examples of these  extensions are the following: PBW extensions defined by Bell and Goodearl (which include the classical commutative polynomial rings, universal enveloping algebra of a Lie algebra, and others); some operator algebras (for example, the algebra of linear partial differential operators, the algebra of linear partial shift operators, the algebra of linear partial difference operators, the algebra of linear partial $q$-dilation operators, and the algebra of linear partial q-differential operators); the class of di\-ffu\-sion algebras; Weyl algebras; additive analogue of the Weyl algebra; multiplicative analogue of the Weyl algebra; some quantum Weyl algebras as $A_2(J_{a,b})$; the quantum algebra $\cU'(\mathfrak{so}(3,\Bbbk))$; the family of 3-dimensional skew polynomial algebras (there are exactly fifteen of these algebras, see \cite{ReyesSuarez2017PBWBasesFarEast}); Dispin algebra $\cU(osp(1,2))$; Woronowicz algebra $\cW_v(\mathfrak{sl}(2,\Bbbk))$; the complex algebra $V_q(\mathfrak{sl}_3(\mathbb{C}))$; $q$-Heisenberg algebra ${\bf H}_n(q)$; the Hayashi algebra $W_q(J)$, and several algebras of quantum physics (for instance, Weyl algebras, additive and multiplicative analogue of the Weyl algebra, quantum Weyl algebras, q-Heisenberg algebra, and others). For a detailed reference of each one of these algebras, see \cite{Reyes2013PhD}, \cite{LezamaReyes2014}, \cite{ReyesSuarez2017Momento}, \cite{Reyes2018} and \cite{ReyesSuarez2018Beit}.
\end{examples}
\begin{remark}
In a forthcoming paper we will present an algorithm to find explicitly the polynomial $f(s,t)$ established in Theorem \ref{moreimportant}.
\end{remark}

\end{document}